\documentstyle[11pt]{article}
\newtheorem{theo}{Theorem}[section]
\newtheorem{prop}[theo]{Proposition}
\newtheorem{lemma}[theo]{Lemma}

\newtheorem{quest}[theo]{Question}

\newtheorem{definition}[theo]{Definition}
\newtheorem{clm}[theo]{Claim}

\newcommand{\ignore}[1]{}

\begin{document}

\title{Distributed Corruption Detection in Networks}

\author{Noga Alon \thanks{Department of Mathematics,
 Princeton University, and 
 Schools of Mathematics and Computer Science,
 Tel Aviv University, 
{\tt nogaa@tau.ac.il}.  Research supported in part by 
NSF grant DMS-1855464, ISF grant 281/17,
BSF grant 2018267
and the Simons Foundation.
}
\and Elchanan Mossel\thanks{Department of Mathematics, 
  Massachusetts Institute for Technology
 {\tt elmos@mit.edu}. Research supported in part by
NSF awards DMS-1737944, ONR N00014-16-1-2227, N00014-17-1-2598, ARO MURI W911NF1910217 and Simons Investigator award (622132). }
 \and 
Robin Pemantle\thanks{Department of Mathematics, 
University of Pennsylvania,
209 South 33rd Street,
Philadelphia, PA 19104, USA, {\tt pemantle@math.upenn.edu}.  
Research supported in part by 
NSF grant \# DMS-1209117.}}

\maketitle

\begin{abstract}
We consider the problem of
distributed  corruption detection in     
networks. In this model
each vertex of a directed graph is     
either truthful or corrupt. 
Each vertex reports
the type (truthful or corrupt) of each of its outneighbors.
If it is truthful,
it reports the truth, whereas if it is corrupt, 
it reports adversarially.
This model, first considered by Preparata, Metze and Chien in 1967,
motivated by the desire to identify
the faulty components of a digital system by having the other components
checking them, became known as the PMC model. The
main known results for this model
characterize networks in which \emph{all} corrupt (that is, faulty)
vertices can be identified,
when there is a known upper bound on their number.

We are interested in
networks in which the identity of a \emph{large fraction}
of the vertices can be identified. 

It is known that in the PMC model,   
in order to identify all corrupt
vertices when their number is $t$, 
all indegrees have to be at least $t$. In contrast, we show
that in $d$ regular-graphs with strong expansion
properties,  
a $1-O(1/d)$ fraction of the corrupt vertices, 
and a $1-O(1/d)$ fraction of the truthful 
vertices can be
identified, whenever there is a majority of truthful vertices.
We also observe that if the graph is very far from being a good expander,
namely, if the deletion of a small set of vertices splits the graph
into small components, then no corruption
detection is possible even if most of the vertices are truthful.
Finally we discuss the algorithmic aspects and the computational
hardness of the problem.
\end{abstract}

\section{Introduction}
We study the problem of
\emph{ corruption detection} in networks. 
Given a network of agents, a subset of whom are corrupt and the rest are truthful, our goal is to 
find as many 
truthful and corrupt   
agents as possible. 
Neighboring vertices audit each other. We assume that 
truthful  agents report the 
type  
of their neighbors accurately. 
We make no assumption on the report of corrupt agents. 
For example, two corrupt neighbors may collude and report each other as
non-corrupt. Similarly, a corrupt vertex may 
prefer to report the status of some of 
its neighbors accurately, hoping that this will establish a 
truthful record for itself. 
We permit that the corrupt agents   
coordinate their actions in an arbitrary fashion. 

The corruption model studied here is identical 
to the model of diagnosable systems 
that was introduced by Perparata, Metze and Chien~\cite{PMC67} 
as a model of a digital system 
with many components that can 
fail. 
It is assumed that components can test 
some other components. The goal in~\cite{PMC67} and
in follow-up  
work,   
including~\cite{HA74,KTA75,KR80},  
is to characterize networks that can detect 
all  
corrupt  
nodes 
as long as their number does not exceed a given value.

Similar models
have been   
introduced and studied in other areas of computer science, 
including Byzantine computing~\cite{LSP} and intrusion 
detection in the security community~\cite{mukherjee1994network}. 
See also the survey~\cite{lgorzata2004survey}.
A byproduct of this line of research is the VLSI chips puzzle discussed
in \cite[Problems~4-6]{CLRS}. In this puzzle there are $n$ supposedly
identical VLSI chips that in principle are capable of
testing each other. A basic test involves two chips,
each chip tests the other and reports whether
it is good or bad. A good chip always
gives an accurate report,   
but the answer of a bad chip cannot be trusted.
The objective is to find a good chip, or all good chips, assuming
more than half of the chips are good, using the minimum possible number
of tests.  This is (an adaptive version of) 
the corruption detection problem
on a complete graph. 

The original motivation for our work
has been  
corruption detection in social
and economic networks,
where the main objective is to understand the
structure of networks that enable one to identify as many of the corrupt
nodes and as many of the truthful ones as possible. 

Examples of such networks may include different
government agencies in a country, the network of banks in the EU or the
network of hospitals in a geographic location.  Our goal is to understand theoretically 
which network structures are more amenable to corruption and which are
more robust against it. In particular, is it possible to design sparse networks that allow to detect corruption even if there are many corrupt nodes? 
Social scientists have studied many 
aspect of corruption networks,
see, e.g., 
\cite{nielsen2003corruption,rock2004comparative,fjeldstad2003fighting}.
However, to the best of our knowledge, prior to this
 work, the effect of the network structure on corruption detection 
has not been systematically studied.     

As noted in follow up work~\cite{JiMoRa:18}, the theoretical models and results discussed here have natural limitations for social and economic networks:
`` Despite our theoretical results, corruption \emph{is} prevalent in many real-world networks, and yet in many scenarios it is not easy to pinpoint even a single truthful node. One reason for that is that some of the 
assumptions do not seem to hold in some real world networks. For example, we assume that audits from the truthful nodes are not only non-malicious, but also perfectly reliable. In practice this assumption is unlikely to be true: many truthful nodes could be non-malicious but simply unable to audit their neighbors accurately. Further assumptions that may not hold in some scenarios include the notion of a central agency that is both uncorrupted and has access to reports from every agency, and possibly even the assumption that the number of corrupt nodes is less than $|V| / 2$". In addition there are many constraints on the network $G$ that may prohibit it to be an expander. 
Despite these shortcomings of our model, our work points to ideal conditions that allow corruption detection. 
While it is perhaps unreasonable to assume that one imposes an expander auditing structure among government agencies, imposing such a structure among more equal entities such as banks, hospitals or universities may be more realistic. Furthermore, our results for directed graphs, suggest such structures even in cases where the auditing relation is not symmetric.

\subsection{Formal Definitions and Main Results} 

\begin{definition}[Digraph] 
For a set $V$, let $V^{(2)}$ denote the set of ordered pairs of distinct elements. 
A {\em digraph} $G=(V,E)$ consists of a set $V$ of nodes (``agents") and a set $E \subset  V^{(2)}$ of directed edges. 
\end{definition} 
\begin{definition}[Type, Report]
Consider a (digraph) $G=(V,E)$ along with a function $\tau : V \to \{\mathsf{c},\mathsf{t}\}$ 
that assigns each node a {\em type}. We call $\tau$ a 
{\em truthfulness assignment} to $G$. 
We call a node $v$ {\em truthful} if $\tau(v) = \mathsf{t}$ and corrupt if $\tau(v) = \mathsf{c}$. 
We write $T = \tau^{-1}(\mathsf{t})$ for the set of truthful nodes and $B =  \tau^{-1}(\mathsf{c})$
for the set of corrupt nodes, so that $V = T \sqcup B$ is a partition of the set of nodes. 
A {\em report} is a function $\rho: E \to \{\mathsf{c},\mathsf{t}\}$. 
A report $\rho: E \to \{\mathsf{c},\mathsf{t}\}$ is {\em compatible} with 
the truthfulness assignment $\tau$ if for each truthful node $u \in T$ and each
directed edge $(u, v) \in E$ we have $\rho(u,v) = \tau(v)$.  
We call $\rho(u,v)$ the type of $v$ reported by $u$. 
We call $\tau$ a {\em valid} truthfulness assignment if $|T| > |B|$. 
We say that a report $\rho$ is {\em feasible} if
it is compatible with at least one valid truthfulness assignment.
 \end{definition}
Since it is common to use the letter $C$ for constants, we prefer to denote the set of corrupt nodes by $B$ (the set of  `` bad" nodes).

The question we address is under what conditions on the digraph
$G$ and the number of truthful vertices it is possible to identify the truthfulness 
status of most nodes 
with certainty. 
It is easy to see
that this is impossible if $|T| \leq |B|$. Indeed, if $V=V_1 \cup
V_2 \cup W$ is a partition of $V$ into $3$ pairwise disjoint sets
where $|V_1|=|V_2|$ (and $W$ may be empty),
then the corrupt agents can ensure that all the reports 
in the two scenarios $T=V_1$, $B=V_2 \cup W$ and $T=V_2$, $B=V_1 \cup W$
will be identical. As there is no common truthful agent  in these
two possibilities, no algorithm can locate a truthful
agent with no error.

Our main result is that if the graph is a good bounded-degree 
directed expander, in the sense described below,
and  we have a majority of truthful agents, then it is possible to
identify the truthfulness status of most of the nodes. 
We recall the following definition of regular spectral expander graphs.
\begin{definition} \label{dec:ndl}
Call a graph an
$(n,d, \lambda)$-graph if it is $d$-regular, has $n$ vertices and all
its eigenvalues besides the top one are in absolute value at most
$\lambda$. 
\end{definition} 
We use the fact that classical construction of expanders 
like the Ramanujan graphs of
\cite{LPS}, \cite{Ma}, are $(n,d,\lambda)$ graphs with $\lambda = 2 \sqrt{d-1}$. 
For the known connection between the 
expansion and pseudo-random properties of graphs and their eigenvalues,
see, e.g., \cite{AC}, \cite{HLW} and the references therein.

First we consider the somewhat simpler case of \emph{undirected graphs}
(the adjacency relation is symmetric).

The main result for the undirected case
is the following.
\begin{theo}
\label{t11}

Let $G = (V,E)$ be a $(n,d,\lambda)$ graph, with $d \geq 100$, which is a Ramanujan graph, i.e. satisfies, $\lambda \leq 2 \sqrt{d-1}$. 
Then, there is an algorithm that given a feasible report $\rho$,  returns subsets $T',B'$ 
such that for every  truth assignment $\tau$ that is compatible with $\rho$, 
 with $T = \tau^{-1}(\mathsf{t})$ and $B =  \tau^{-1}(\mathsf{c})$, and $|T| > |B|$, 
it holds that 
\begin{equation} \label{eq:main_cond_ram}
T' \subset T, \quad B' \subset B, \quad |T \setminus T'| \leq \frac{32}{d} |B| , \quad 
|B \setminus B'| \leq  \frac{32}{d} |B|. 
\end{equation}
Moreover, if there exists a truth assignment compatible with 
$\rho$ with $|T| >  (1 +  \frac{12}{d}) \frac{n}{2} $, then there is a linear-time
algorithm that  identifies subsets $T'$ and 
$B'$ satisfying~(\ref{eq:main_cond_ram}).
\end{theo}

We remark that 
\begin{enumerate}
\item Given any graph $G$ with all degrees bounded by $d$, with at least $2d+1$ corrupt nodes, then there is a set of  $\lfloor |B|/(2d+1) \rfloor$ corrupt nodes and a set of $\lfloor |B|/(2d+1) \rfloor$ truthful nodes who cannot be identified, 
 even if in addition to the report $\rho$, we are given the number of corrupt nodes. Thus, the fractions of nodes we identify is tight up to a constant factor. 
To see that this is the case, consider the following selection of $b$ corrupt nodes:  
choose an arbitrary vertex $v_1$ in the graph and set all its neighbors $N(v_1)$ to be corrupt, then choose $v_2 \notin \{ v_1 \} \cup N(v_1)$ and set all of its neighbors to be corrupt, then choose one of $v_1,v_2$ to be corrupt 
and the other to be truthful. Continue in a similar fashion, defining $N(v_3),N(v_4)$ to be corrupt and choosing one of $v_3,v_4$ to be corrupt and the other truthful. 
Continue until the number of corrupt nodes left, $b'$ satisfies $b' < 2d+1$. 
Set $b'$ of the remaining nodes to be corrupt. 
 Let $r(v,u) = \mathsf{c}$, whenever $v$ is corrupt.  It is then clear that it is impossible to decide which of the $v_i$ is corrupt and which is truthful.
 
\item In the case where $|B|/d = o(n)$,  the theorem allows to recover 
$1-o(1)$ fraction of the good nodes. 

\item 
The algorithm in the proof of the theorem is an
exponential time algorithm if 
we only assume that $|T|>|B|$ (or if we assume that $|T|>(1/2+\mu)n$ for a
very small fixed $\mu=\mu(\lambda,d)$). 

\end{enumerate}

The fact that the detection algorithm is not efficient
when we only assume that $T$ is just a little bit bigger than $B$
is not a coincidence. Indeed, the algorithm described
in the proof of the theorem, presented
in the next sections, provides a set $T$ of more than $n/2$ truthful
agents, which is consistent with the reports obtained, when such a
set exists. We show that the problem of producing
such a set when it exists is $NP$-hard,
even when restricted to bounded-degree expanders (and even if 
we ensure that there is such a set of size at least
$n/2+\eta n$.)
\begin{theo}
\label{t12}
There exist $c > 0, \eta > 0$ 
such that the following promise problem is 
$NP$-hard. The input is a graph $G=(V,E)$ with $|V| = n$, 
which is an $(n,d,c\sqrt{d})$-graph  
along with a report $\rho$. 
The promise is that either
\begin{itemize}
\item 
$\rho$ is compatible with $\tau$ satisfying  
$|T| = |\tau^{-1}(\mathsf{t})| \geq n/2 + \eta n$, or
\item 
Any $\tau$ which is compatible with $\rho$ satisfies 
 $|T| = \tau^{-1}(\mathsf{t})  \leq n/2 - \eta n$. 
\end{itemize}
The objective is to distinguish between the two options above.
\end{theo}
The proof is presented in subsection~\ref{subsec:hardness}.

We also establish in subsection~\ref{subsec:bad_expanders} the following simple statement, 
which shows that at least some
(weak) form of expansion is needed for solving the
corruption detection problem.
\begin{prop}
\label{p13}
Let $G=(V,E)$ be a graph on $n$ vertices so that it is possible to
remove at most $\epsilon n$ vertices of $G$  and get a graph in which 
each   
connected component is of size at most $\epsilon n$. Then given a report 
$\rho$ compatible with an assignment $\tau$, with $T = \tau^{-1}(\mathsf{t})$ and $B =  \tau^{-1}(\mathsf{c})$,
and $|T| \geq (1-2 \epsilon) n$ it is impossible to
identify even a single member $t \in T$ from the reports of
all vertices. In particular, this is the case for planar graphs or
graphs with a fixed excluded minor even if
$\epsilon=\Theta(n^{-1/3})$.
\end{prop}

Note that there is still a significant gap between the expansion properties
that suffice for solving the detection problem, described in Theorem~\ref{t11},  
and the conditions in the last proposition that are
necessary for such a solution. It will be interesting to obtain tighter
relations between expansion and corruption detection. This is
further discussed in section 4.

\subsection{Results for directed graphs}
In this subsection we consider directed graphs (digraphs). 
This is motivated by the fact
that in various auditing situations it is 
unnatural to allow $u$ to inspect $v$ whenever 
$v$ inspects $u$. In fact, it may even be desirable not to allow any short cycles in the directed inspection graph. 

For a set of vertices $A$, in a graph $G$, define:
\begin{equation} \label{def:na}
N(A) := \{ v : \exists u \in A, \{u,v\} \in E \},
\end{equation}

For a digraph $G=(V,E)$, let 
\[
N^{+}(U) = \{ (v : \exists u \in U, (u,v) \in E\}, \quad 
N^{-}(U) = \{ (v : \exists u \in U, (v,u) \in E\}.
\]
Note that $N(U), N^{+}(U)$ and $N^{-}(U)$ may contain elements from $U$. 
We first present an explicit construction of regular directed 
expanders which expand at three 
different scales as in the undirected case. 

\begin{prop}
\label{p14} 
There exists $c_1,c_2 > 0$ and  infinitely many values of 
$d$ for which there are infinitely many values of $n$ and an 
explicit construction of oriented $3(d-6)$-regular
graphs $G=([n],E)$ which satisfy 
\begin{itemize}
\item The girth of $G$ as an undirected graph is at least $2 \log_{d}(n) / 9$. 
\item If $A \subset [n]$ is of size at most $n/(9 d)$ then $|N^{+}(A)| > c_1 d |A|$ and $|N^{-}(A)| > c_1 d |A|$. 
\item For any set $A \subset [n]$ of size at most $n/4$, it holds that $|N^{+}(A)| > 2 |A|$ and $|N^{-}(A)| > 2|A|$. 
\item If $A,B \subset [n]$ with $|A| \geq c_2 n /d$, and $|B| \geq n/4$, then there is a directed edge from $A$ to $B$ and a directed edge from $B$ to $A$. 
\end{itemize} 
\end{prop} 

The proof, presented in Section~\ref{sec:directed}, is based 
on ``packing" three different group based expanders 
to obtain the expansion in the different scales.  Given 
the directed expanders constructed in Proposition~\ref{p14}, we prove the following directed analogue of Theorem~\ref{t11}.

\begin{theo}
\label{t15}
Consider the graphs constructed in Proposition~\ref{p14} . 
There exist constant $c_3(c_1,c_2),c_4(c_1,c_2)$, for which the following holds. 

There is an algorithm that given a feasible report $\rho$,  returns subsets $T',B'$ 
such that for every  truth assignment $\tau$ that is compatible with $\rho$, 
 with $T = \tau^{-1}(\mathsf{t})$ and $B =  \tau^{-1}(\mathsf{c})$, and $|T| > |B|$, 
it holds that 
\begin{equation} \label{eq:main_cond_dir}
T' \subset T, \quad B' \subset B, \quad |T \setminus T'| \leq \frac{c_3}{d} |B| , \quad 
|B \setminus B'| \leq  \frac{c_3}{d} |B|. 
\end{equation}
Moreover, if there exists a truth assignment compatible with 
$\rho$ with $|T| >  (1 +  \frac{c_4}{d}) \frac{n}{2} $, then there is a linear-time
algorithm that  identifies subsets $T' \subset T$ and a subset
$B' \subset B$ satisfying~(\ref{eq:main_cond_dir}).
If we only assume that $|T|>|B|$ 
then the detection algorithm is exponential.
\end{theo}

\subsection{Relation to Distributed Computing} 
The model presented in the paper requires a central agency that obtains the report $\rho$ from all nodes. This is in contrast to models in 
distributed computing and byzantine agreement where all nodes 
only communicate with neighboring nodes. We note however that the
following holds.

\begin{prop} \label{p:distributed}
In the setting of  of Theorem~\ref{t11} when 
assuming each node has a unique identity that is known to all of its neighbors, 
there is a distributed algorithm on the graph $G$, where all nodes in $T'$ 
output the sets $T'$ and $B'$ in~(\ref{eq:main_cond}).
\end{prop} 

\subsection{Related Work} 

The vast literature on corruption detection in computer science,
and in particular on the diagnosable system problem and the PMC model
introduced in \cite{PMC67},
deals  
either with the problem of identifying all
corrupt nodes, or with that of identifying a single corrupt node. As observed
in \cite{PMC67},  a necessary condition 
for the identification of \emph{all} corrupt nodes in a network with $t$ 
corrupt nodes is that the \emph{minimal indegree} 
in the network is at least $t$.
Therefore, if the number of corrupt nodes is linear in the total number
of vertices, all indegrees have to be linear, and the total
number of edges has to be quadratic.

The main contribution of the present paper is 
a proof that the number of required edges may be much smaller
when relaxing the requirement of 
identifying \emph{all} corrupt nodes and replacing it by the requirement of the 
identification of the truth status of most of the nodes.
By relaxing the requirement as above 
we are able to study bounded-degree graphs, in particular $d$-regular graphs. 
Our main new result is that a linear number
of edges ensures the recovery of the the truth status of a large fraction of the nodes, 
provided the graph is a sufficiently strong expander.
It was shown already
in \cite{PMC67} that   a linear number of edges suffices
to ensure the detection
of \emph{a single} corrupt vertex. We show that such a small number of
edges suffices to determine the types of a $1-O(1/d)$ fraction of the vertices, even when
the number of truthful vertices exceeds that of corrupt ones by only
$1$.

In the context of Byzantine agreement it was discovered already
in the 80s that expanders allow ``almost everywhere agreement''. This
was first established by Dwork et. al. in \cite{DPPU} and  further 
developed in subsequent work, see, in particular \cite{GO} and its 
references. In \cite{DPPU} it is shown that one can achieve a relaxed notion of 
Byzantine agreement, where a large fraction of non-faulty nodes agree on a value for random regular graphs. Like in our results, in the bounded degree case, this fraction is bounded away from $1$ unless the number of corrupt nodes is sublinear. 
Moreover, the results of \cite{DPPU} require that the number of 
faulty nodes is a small fraction 
(much smaller than $1/2$) of the total number of nodes. 

It is therefore not very surprising that graph 
expansion is relevant to corruption detection. It is, however, interesting
to note that in sufficiently strong expanders it is possible to identify
most of the truthful and most of the corrupt agents even if the number of
truthful agents exceeds the number of corrupt ones by only $1$.

\subsection{Techniques} 

The proof of Theorem~\ref{t11} rely crucially on the fact that $(n,d,O(\sqrt{d}))$-graphs expand well at different scales.  
\begin{itemize}
\item Every set $A$ of size $\Omega(n/d)$ and every set $B$ of size $n/4$ have at least one edge between them, 
\item Every set $A$ of size $O(n/d)$ has $|N(A) \setminus A| \geq |A|$.  
\item Every set such that $|A \cup N(A)| \leq 3n/4$ satisfies $|N(A)| \geq \Omega(d) |A|$. 
\end{itemize} 

For the directed case, we prove the existence of directed expanders of high girth with analogous properties in 
Proposition~\ref{p14}.

We observe that a certain weak expansion property, 
namely, the nonexistence of small separators, is 
necessary for corruption detection.
Combining this observation with the 
planar separator theorem of Lipton and Tarjan and its extensions
we conclude that planar graphs and graphs with a fixed excluded minor 
are not good for corruption detection. 

Finally we discuss the algorithmic aspects of our problem 
using results about hardness of approximation.

\section{Proofs} \label{sec:proofs} 

In this section we present the proofs of all results besides the ones
for directed graphs. These appear in the next section.

\subsection{Undirected graphs} \label{subsec:undirected}
We first state the following generalization of Theorem~\ref{t11}. 
Theorem~\ref{t11} follows when considering graphs 
with $\lambda \leq 2 \sqrt{d-1}  ~(~ < 2 \sqrt d)$. 
\begin{theo}
\label{t21}

Let $G = (V,E)$ be a $(n,d,\lambda)$ graph, where
$d^2 \geq 24 \lambda^2$.  
Then, given a report $\rho$ which is compatible with 
a truth assignment $\tau$ with $T = \tau^{-1}(\mathsf{t})$ and $B =  \tau^{-1}(\mathsf{c})$, and $|T| > |B|$, 
we can identify  a subset $T' \subset T$ and a subset
$B' \subset B$
so that 
\begin{equation} \label{eq:main_cond}
|T \setminus T'| \leq \frac{8 \lambda^2}{d^2} |B| , \quad 
|B \setminus B'| \leq  \frac{8 \lambda^2}{d^2} |B|.
\end{equation}

Moreover, if $|T|>(1/2+3\lambda^2/d^2)n$ then there is a linear-time
algorithm that  identifies a subset $T' \subset T$ and a subset
$B' \subset B$ satisfying~(\ref{eq:main_cond}).

\end{theo}

For a positive $\delta<1/8$
call a graph $G=(V,E)$ on a set of $n$ vertices a
$\delta$-{\em good}
expander if any set $U$ of at most $2 \delta n$ vertices has 
more than $|U|$ neighbors outside $U$, and there is an edge between
any pair of sets of vertices provided one of them is of
size at least $\delta n$ and the other is of size at least 
$n/4$.
We recall how standard results about expanders 
imply that $(n,d,O(\sqrt{d}))$-graphs are $\delta$ good 
for $\delta = \Theta(1/d)$. 

\begin{lemma}[\cite{ABNNR}, Corollary~1] \label{cor:expand}
Let $G=(V,E)$ be an $(n,d,\lambda)$-graph and let $B \subset V$ be of size 
$b n$. Let $t n$ denote the size of the set $V \setminus N(B)$, of nodes that have no neighbors in $B$. Then:
\[
t \leq \frac{\lambda^2}{d^2} \frac{1-b}{b}.
\]
\end{lemma}


\begin{lemma} \label{lem:good_expanders} 
Any $(n,d,\lambda)$-graph
in which $\frac{3 \lambda^2}{d^2} \leq \delta \leq 1/8$, is a 
$\delta$-good expander.
\end{lemma} 

{\bf Proof:} 
Let $U$ be of size at least $\delta n$ and suppose that that there are no edges between 
$U$ and set $B$ of size $n/4$. Then if $u = |U|/n$, by 
Lemma~\ref{cor:expand} it follows that 
\[
u \leq \frac{3 \lambda^2}{d^2},
\]
which is a contradiction. 

Similarly, let $U$ be a set of size  $ u n \leq 2 \delta n$ and suppose that 
$|N(U) \setminus U| < u n$. Then the set $B = V \setminus 
(U \cup N(U))$ is of size 
at least $b n$ where $b > 1- 2 u \geq 1/2$. 
This implies by Lemma~\ref{cor:expand} that
\[
u \leq \frac{\lambda^2}{d^2} \frac{1-b}{b} \leq \frac{\delta}{3} 4 u \leq u/2,
\]
which is a contradiction.
\hfill $\Box$

The other property of expanders we will need is that small sets expand by a factor of 
$\Omega(d)$. In particular:

\begin{lemma} \label{lem:small_sets} 
Let $G=(V,E)$ be an $(n,d,\lambda)$-graph. 
Let $A \subset V$ be such that $|A \cup N(A)| \leq 3n/4$.
Then 
\begin{equation} \label{eq:small_sets}
| A \cup N(A)| \geq \frac{d^2}{4 \lambda^2} |A|.
\end{equation}
\end{lemma} 

{\bf Proof:} 
The proof is similar. Let $|A| = a n$. 
Let $B = V \setminus (A \cup N(A))$ and put 
$|B| = (1- c a) n \geq n/4$. Since there are no edges between $A$ and $B$, 
it follows that 
\[
a \leq \frac{\lambda^2}{d^2}  \frac{c a}{1 - c a} \leq 4 c a \frac{\lambda^2}{d^2},
\]
and so 
\[
c \geq  \frac{d^2}{4 \lambda^2},
\]
which is the same as~(\ref{eq:small_sets})
\hfill $\Box$

\noindent
{\bf Proof of \ref{Theorem}{t21}:}\,
Let $G$ be an $(n,d,\lambda)$-graph, where
$d^2 \geq 24 \lambda^2$.

Let $\tau$ be a truth assignment consistent with $\rho$ and let $T = \tau^{-1}(\mathsf{t})$ and 
$B = \tau^{-1}(\mathsf{c})$. 
Note that by lemma~\ref{lem:good_expanders}, $G$ is $\delta$ good, where 
$\delta = 3\lambda^2/d^2$. 

Let $H$ be the spanning  subgraph of $G$ in which $(u,v) \in E$ iff
$\rho(u,v) = \rho(v,u) = \mathsf{t}$. 
Let $V_1,V_2, \ldots ,V_s$ 
be the sets of vertices of the connected  components of 
$H$. 
\begin{clm}
\label{c21}
All the vertices of each $V_i$ are of the same type, that is, for
each $1 \leq i \leq s$, either $V_i \subset T$ or $V_i \subset B$.
\end{clm}

\noindent
{\bf Proof:}\, Suppose $u$ and $v$ are neighbors in $H$. If $u \in T$
then $v \in T$ (as $u$ reports so). If $u \in B$, then 
$v \in B$ (as $v$ reports that $u \in T$).  \hfill $\Box$

Call a component of $H$ truthful if it is a subset of $T$, 
otherwise   
it is a subset of $B$ and we call it corrupt. 

Let $H'$ be the induced subgraph of $G$ on the set $T$ of all 
truthful  
vertices.
\begin{clm}
\label{c22}
Any connected component of $H'$ is also a connected component of 
$H$.
\end{clm}

\noindent
{\bf Proof:}\, If $u,v \in T$ are adjacent in $G$ (and hence in
$H'$), they are
adjacent in $H$ as well, by definition and by the fact that each of
them reports honestly about its neighbors. Thus each component $C'$ of
$H'$ is contained in a component $C$ of $H$. However, no
$v \in T$ is adjacent in $H$ to a vertex $w \in B$, implying
that in fact $C'=C$ and establishing the
assertion of the claim. \hfill $\Box$
\begin{clm}
\label{c23}
The graph $H'$ contains a connected component of size at least
$|T|-\delta n > (1/2-\delta)n$.
\end{clm}

\noindent
{\bf Proof:}\, Assume this is false and the largest connected
component of $H'$ is on a set of vertices $U_1$ of size smaller
than $|T|-\delta n$. Since the total number
of vertices of $H'$ is $|T|>n/2$, it is easy to check that one can
split the connected components of $H'$ into two disjoint sets,
each of total size  at least $\delta n$. However, the bigger among the two
is of size bigger than $n/4$, and hence, since $G$ is a
$\delta$-good expander, there is an edge of $G$ between  the two
groups. This is impossible, as it means that there is an edge
of $G$ between two distinct connected
components of $H'$. \hfill $\Box$

The analysis so far allow us to prove the easy part of the theorem.
\begin{clm}
If  $|T|>(1/2+\delta)n$ then there exists a linear-time algorithm 
which finds $T' \subset T$ and $B' \subset B$ such that 
\[
|T \setminus T'| \leq 8 \frac{\lambda^2}{d^2} |B| , \quad 
|B \setminus B'| \leq  8 \frac{\lambda^2}{d^2} |B|. 
\]
\end{clm}
\vspace{0.1cm}

\noindent
{\bf Proof}:\,
Note that if $|T|>(1/2+\delta)n$ then Claim \ref{c23} implies that 
$H$ must contain a connected component of size bigger than
$n/2$, which must be truthful. Denote the vertices of this component by $T'$ and 
$B' = N(T')$ and $R = V \setminus (B' \cup T')$. Clearly $B' \subset B$ and 
by Lemma~\ref{lem:small_sets} it follows that
\[
|R \cup B'| \geq \frac{d^2}{4 \lambda^2} |R|
\]
 which by the assumption that $d^2 \geq 24 \lambda^2$ that 
gives $|R|+|B'| \geq 6|R|$ and hence $|R \cup B'| \leq \frac{6}{5}|B'|$
implies 
 \[
|T \setminus T'| \leq |R| <
5 \frac{\lambda^2}{d^2} |B| , \quad |B \setminus B'| \leq 
|R|<  5 \frac{\lambda^2}{d^2} |B|, 
\]
establishing the required inequality (with room to spare).

Thus, if this is the case, the set $T'$ is found  
by the simple, linear-time algorithm that computes the connected
components of $H$. Furthermore, the set $B'$ can also be found by computing the vertex boundary of $H$ which is easily computed in linear time as well. 
 \hfill $\Box$

It remains to show that even if we only
assume that $|T|>n/2$ then we can still identify correctly most of the
truthful and corrupt vertices. 
We proceed with the proof of this stronger
statement.

By Claim \ref{c22} and Claim \ref{c23} it follows that $H$ contains 
at least one
connected component of size at least $(1/2-\delta)n \geq 3/8 n$.
If $H$ contains only one such component, then this component  
must consist of truthful agents, and we can identify all of them.
Having these truthful vertices, we also
know the types of all their neighbors. By the assumption on $G$
this gives the types of all vertices but less than $\delta n$, thus establishing
(\ref{eq:main_cond}).

Otherwise, there is another connected component of size at least
$(1/2-\delta)n$, and as there is no room for more than two such
components, there are exactly two of them, say $V_1$ and $V_2$. Note
that by the expansion properties of $G$ there are edges of $G$ between $V_1$
and $V_2$, and hence it is impossible that both of them are
truthful components.  As one of them must be truthful, it follows
that exactly one of $V_1$ and $V_2$ is a truthful component and the 
other corrupt. We next show that we can identify the types of both
components.

Construct an auxiliary weighted graph $S$ on the set of vertices 
$1,2,
\ldots ,s$ representing the connected components  $V_1, V_2
\ldots ,V_s$ as follows. The weight $w_i$ of $i$ is defined
by $w_i=\frac{|V_i|}{|V|}$. Two vertices $i$ and $j$ are
connected iff there is at least one edge  of $G$ that connects a
vertex in $V_i$  with one in $V_j$. Call an independent set in
the graph $S$ \emph{large} if its total weight is bigger than $1/2$. Note
that by the discussion above $T$ must be a union of the form
$T=\bigcup_{i \in I} V_i$,  
where $I$ is a large independent set in the
graph $S$. In order to complete the argument we prove the following.
\begin{clm}
\label{c24}
Either there is no 
large independent set  in $S$ containing $1$, or there is no
large independent set in $S$ containing $2$.
\end{clm}
\vspace{0.1cm}

\noindent
{\bf Proof:}\,  Assume this is false, and let $I_1$ be a large 
independent set in $S$ containing  $1$, and $I_2$ a large 
independent set in $S$ containing $2$. To get a contradiction 
we show that for
$ w(I_1)=\sum_{i \in I_1}w_i$ and $w(I_2) =\sum_{i \in I_2} w_i$
we have $w(I_1) + w(I_2)  \leq 1$ (and hence it is impossible that
each of them has total weight bigger than a half).

To prove the above note, first, that the two vertices $1$ and $2$
of $S$ are connected (as each corresponds to a set of more than
$(1/2-\delta)n$  vertices of $G$, hence there are edges of $G$ 
connecting $V_1$ and $V_2$). Therefore $I_1$ must contain 
$1$ but not $2$, and $I_2$ contains $2$ but not $1$.

If there are any vertices $i$ of $S$ connected in $S$ both to 
$1$ and to $2$, then these vertices belong to neither $I_1$ nor $I_2$,
as these are independent sets. Similarly, if a vertex $i$ is
connected to $1$ but not to $2$, then it can belong to
$I_2$ but not to $I_1$, and the symmetric statement holds for
vertices connected to $2$  but not to $1$. So far we have discussed
only vertices that can belong to at most one of the two independent
sets $I_1$ and $I_2$. If this is the case for all the vertices of
$S$, then each of them contributes its weight only to one of the
two sets and their total weight would thus be at most $1$, implying
that it cannot be that the weight of each of them is bigger than
$1/2$, and completing the proof of the claim. It thus remains to
deal with the vertices of $S$ that  belong to both $I_1$ and $I_2$.
Let $J \subset \{3,4, \ldots  ,s\}$ be the set of all these vertices.
Note, first, that the total weight of the vertices in $J$ is at
most $2 \delta$, as the total weight of $1$ and $2$ is at least
$2 (1/2-\delta)=1-2 \delta$.  Note also that by the discussion
above each $j \in J$ is not a neighbor of $1$ or of $2$. 
By the assumption about the expander $G$ the total weight of the
vertices that are neighbors of vertices in $J$ and do not belong to
$J$ is bigger than the total weight of the vertices in $J$.
Indeed, this is the case as the number of neighbors in $G$
of the set $\cup_{j \in J} V_j$ that do not lie in this set
is bigger than the size of the set. We thus conclude that
if $J'=N_S(J)-J$ denotes the set of neighbors of $J$ that do not belong
to $J$, then the total weight of the vertices in $J'$ exceeds the
total weight of the vertices in $J$, and the vertices in $J'$ 
belong to neither $I_1$ nor $I_2$.  We have thus proved that the
sum of weights of the two independent sets $I_1$ and $I_2$
satisfies
$$
w(I_1)+w(I_2) \leq 2 w(J) +(1-w(J)-w(J')) \leq w(J)+w(J')
+(1-w(J)-w(J')) =1
$$
contradicting the fact that both $I_1$ and $I_2$ are large. This
completes the proof of the claim.  \hfill $\Box$

By \ref{Claim}{c24} we conclude that one can identify the types of
the components $V_1$ and $V_2$. This means that we can identify
at least $(1/2-\delta)n$ truthful vertices with no error. Recall
that this is
the case also when $H$ has only one connected component of size
at least $(1/2-\delta)n$. Having these truthful vertices, we also
know the types of all their neighbors. By the assumption on $G$
this gives the types of all vertices but less than $\delta n$,
completing the proof of (\ref{eq:main_cond}).   
\vspace{0.1cm}

\noindent
As noted above the algorithm described in the proof above
is a linear time algorithm provided $|T|>(1/2+\delta)n$.
However, if we only assume that $|T|>|B|$ the proof
provides only a non-efficient algorithm  for deciding the types of
the components $V_1$ and $V_2$. Indeed, we have to compute the
maximum weight of an independent set containing $1$
in the weighted graph $S$,
and the maximum weight of an independent set 
containing $2$. By the proof above, exactly one of this maxima 
is larger than $1/2$, providing the required types. 
\hfill $\Box$

The proof of Proposition~\ref{p:distributed} follows from the fact that it deals with the case of a large component of truthful nodes who can communicate with each other. 

\noindent
{\bf Proof:}
For simplicity, consider a synchronized protocol. In the protocol, in each round, each node transmits to all of its neighbors, the identity of all nodes it recognizes as truthful and as corrupt.  
A truthful node will add to its record 
the identity of all nodes it receives from truthful neighbors. 
Recall that there is exactly one component $T'$ of size $\geq (1/2-\delta)n$ of truthful vertices, that $T'$ is the set of 
vertices in this component and $B'$ is $N(T')$. It is clear that all element in $T'$ who will follow the protocol, will recognize all elements of $T'$ as truthful and all elements of $B'$ as corrupt. 
\hfill $\Box$

\subsection{Hardness} \label{subsec:hardness}

In this subsection we prove Theorem \ref{t12} which explains 
the non-efficiency of the algorithm in the proof of Theorem \ref{t11}.
\vspace{0.2cm}

\noindent
{\bf Proof of Theorem \ref{t12}:}\,
The proof is based on the following result~\cite{BK99}: there exist
constants $b<a<1/2$ such that deciding if a graph $H$ on $m$ vertices,
all of whose degrees are bounded by $4$, has a maximum independent set
of size at least $(a+b) m$ or at most $(a-b) m$ is $NP$-hard. It is easy to see that in fact one can assume that 
$H$ is $4$-regular. 
 
Let $G'$ be an $(n,d-4,c_1\sqrt{d})$-graph with vertex set $V$, where $c \sqrt{d} \geq 8$.  
Split the vertices into $3$ disjoint sets $V_1,V_2,V_3$ of size $\Omega(n)$, 
where $V_3$ is an independent set in $G'$ of size $m$, all its neighbors are in $V_2$, $|V_1|=n/2-am$ and $|V_2|=n/2-m+am$, 
where $m = \Omega(n)$. Write $\eta$ for the constant satisfying $b m=\eta n$. 
Add on $V_3$ a bounded-degree graph $H$ as above, in which it is hard to
decide if the maximum independent set is of size at least $(a+b)m$ or at
most $(a-b)m$. That is, identify the set of vertices of $H$ with $V_3$
and add edges between the vertices of $V_3$ as in $H$. 
Add a $4$-regular graph on $V_1$ and another one on $V_2$. 
Call the resulting]
graph $G$ and note that it is an $(n,d,c_1\sqrt{d} + 4)=
(n,d,c\sqrt d)$ graph. 

The reports of the vertices are as follows.  Each vertex in $V_1$
reports true on each neighbor it has in $V_1$, and corrupt on  any other
neighbor. Similarly, each vertex of $V_2$ reports true on any neighbor
it has in $V_2$ and corrupt on any other neighbor, and each vertex in $V_3$
reports corrupt on all its neighbors.  Note that with these reports the
connected components of the graph $H$ in the proof of Theorem \ref{t11}
are $V_1$, $V_2$ and every singleton in $V_3$.

It is easy to check that here if $H$ has an independent set $I$ of size
at least $(a+b)m$, then $G$ has a set $T$ of truthful  vertices of size
at least $n/2+bm$, namely, the set $I \cup V_1$, which is consistent with
all reports. If $H$ has no independent set of size bigger than $(a-b)m$,
then $G$ does not admit any set $T$ of truthful vertices of size bigger
than $n/2 - bm$ consistent with all reports. This completes the proof.
\hfill $\Box$
\subsection{Graphs with small separators} \label{subsec:bad_expanders}

In this subsection we describe the simple proof of Proposition \ref{p13}. 
\vspace{0.3cm}

\noindent
{\bf Proof of Proposition \ref{p13}:}\,
Let $B'$ be a set of at most $\epsilon n$ vertices of $G$ whose
removal splits $G$ into connected components with vertex classes
$V_1,V_2, \ldots ,V_s$, each of size at most $\epsilon n$. Consider
the following $s$ possible scenarios $R_i$, for $1 \leq i \leq s$.
\vspace{0.1cm}

\noindent
$R_i$:  the set of corrupt vertices is $B=B' \cup V_i$, all the others
are 
truthful   
vertices. The vertices in $B'$ report that all their
neighbors are corrupt. The vertices in $V_i$ report that their
neighbors in $V_i$ are in $T$, and that all their other neighbors 
are in $B$. (The truthful vertices, of course, report truthfully
about all their neighbors).
\vspace{0.1cm}

\noindent
It is not difficult to check that in all these $s$ scenarios, all
vertices make exactly the same reports. On the other hand, there is
no vertex of $G$ that is truthful in all these scenarios, hence 
it is impossible to identify a truthful vertex with no
error. Since the number of corrupt vertices in all scenarios is at most
$2 \epsilon n$, the first assertion of the theorem follows. The
claim regarding planar graphs and graphs with excluded minors
follows from the results in \cite{LT}, \cite{AST}.
\hfill $\Box$

\section{Directed Graphs} \label{sec:directed}

\subsection{Construction of Directed Expanders}

Here we provide the proofs for the case of directed graphs. We start with
the proof of existence of explicit directed expanders that expand 
in three scales. 
We will use the following lemma for undirected graphs from~\cite{AKRS}. 
\begin{lemma} \label{lem:frac_expand}[Lemma 3.6~\cite{AKRS}]
Let $G=(V,E)$ be an $(n,d,\lambda)$ graph. Let $6/d < \gamma < 1$. Let $A,X \subset V$ be sets of vertices such that 
\begin{itemize}
\item $|A| \leq \frac{\gamma}{2(d+1)} n$,
\item for all $v \in A$, it holds that $|\{ w \in X : (w,v) \in E\}| \geq \gamma(d+1)$.
\end{itemize}
Then: 
\[
|X| \geq \frac{\gamma^2 d^2}{9 \lambda^2} |A|.
\]
\end{lemma} 

We will also use the following variant of Lemma~\ref{lem:good_expanders} whose proof is similar. 

\begin{lemma} \label{lem:good_expanders2} 
Any $(n,d,\lambda)$-graph
in which $16 \frac{\lambda^2}{d^2} \leq \delta$ satisfies that 
any set of size at least $\delta n / 2$ and any set of size $n/8$ have at least one edge between them.
\end{lemma} 

{\bf Proof:}
Let $U$ be of size at least $\delta n/2$ and suppose that that there are no edges between 
$U$ and set $B$ of size $n/8$. Then if $u = |U|/n$, by Corollary~\ref{cor:expand} it follows that 
\[
u \leq 7 \frac{\lambda^2}{d^2},
\]
which is a contradiction. 

\hfill $\Box$

\noindent
{\bf Proof of Proposition \ref{p14}:}\,
Let $G'=([n],E')$ be a $d$-regular 
undirected non-bipartite Ramanujan Cayley graph as constructed in
\cite{LPS} or \cite{Ma}. 
This is a Cayley graph of a finite group $\Gamma$ of size $n$, with respect to a set $S'$ of $d$ generators  
$S' = \{a_1, a_1^{-1}, a_2, a_2^{-1},\ldots, a_{d/2}, a_{d/2}^{-1}\}$. 
This graph is known to have girth at least $2 \log_d(n)/3 $, which is the same as saying there is no nontrivial word of the generators of length less than $2 \log_d(n)/3$ that is the identity. 

Let $T = \{a_4,a_4^{-1},\ldots,a_{d/2},a_{d/2}^{-1}\}$. 
Note that the Cayley graph of $\Gamma$ with respect to $T$ is a $d-6$ regular graph whose second eigenvalue is at most $2 \sqrt{d} + 6$. Thus if $ d \geq 36$ this is an 
$(n,d-6,3 \sqrt{d})$-graph. 

Now for $1 \leq i \leq 3$, let $T_i = a_i^{-1} T a_i$. 
Then clearly, $G_i$, the Cayley graph  of $\Gamma$ with respect to $T_i$ is also an $(n,d-6,3 \sqrt{d})$-graph. Moreover, if $S = T_1 \cup T_2 \cup T_3$, then the Cayley graph $H$ of 
$\Gamma$ with respect to $S$ has girth at least $\frac{2}{9} \frac{\log n}{\log d}$, since a nontrivial word in $S$ of length $k$ corresponds to a non-trivial word of length at most $3 k$ in $S'$. 
Note also that $G'$ is $3(d-6)$-regular. 

The desired graph, $G=([n],E)$, will be obtained by assigning directions to the edges $\bar{E}$ of $H$ as follows: 
\begin{itemize}
\item If $\{a,b\}$ is an edge of $G_1$ then orient it from $a$ to $b$ if $a > b$ and from $b$ to $a$ if $b > a$.
\item If $\{a,b\}$ is an edge of $G_2$ then orient it from $b$ to $a$ if $a > b$ and from $a$ to $b$ if $b > a$.
\item In the graph $G_3$ all the degrees are even. Pick an orientation of the edges of 
$G_3$ by picking a directed Eulerian cycle and orienting the edges according to it. In particular, note that every in-degree and out-degree is 
exactly $d/2-3$ in $G_3$. 
\end{itemize} 

We now verify the expansion at the three different scales.
First, we apply Lemma~\ref{lem:frac_expand} to the graph $G_3$ and sets $A$ of size $|A| \leq n/(9d)$, 
 and $\gamma = 1/3$ and obtain that 
\[
|N^{+}(A)| \geq \frac{(d-6)^2}{81 (3 {\sqrt d})^2} |A| \geq \frac{d}{400} |A|,
\]
for $d \geq 36$. The proof for $N^{-}(A)$ is identical.

Next, let $A$ and $B$ be two sets of vertices with  $n/16 \geq |A| \geq 200 n / d$ and $|B| \geq n/4$. Let $m(A)$ and $m(B)$ denote the medians of the sets of numbers given by $A$ and $B$.  Then either $m(A) \geq m(B)$ or 
$m(B) \geq m(A)$. If  $m(A) \geq m(B)$, let $A' = \{ v \in A : v \geq m(A)\}$ and $B' = \{v \in B : v \leq m(B)\}$.
Then $|A'| \geq 100 n / d$ and $|B'| \geq n/8$ and all elements of $A'$ are bigger than all elements of $B'$. 
Thus by Lemma~\ref{lem:good_expanders2} applied to $G_1$, with $\delta = 200/d$, 
$\lambda = 3 \sqrt{d}$ and $d-6 \geq 30$, 
there is at least one edge in 
 $G_1$ connecting $A'$ and $B'$. 
This is a directed edge from $A$ to $B$. A symmetric argument applies if $m(A) \leq m(B)$. 

Note that the last statement implies that if $n/4 \geq |A| \geq 200 n / d$ then the set of non-neighbors of $A$ is of size at most $n/4$ and therefore $|N(A)| \geq 3 n/4 > 2 |A|$. 

For sets $A$ with $n / (10d) \leq |A| \leq 200 n / d$, we have 
\[
|N(A)| \geq  \frac{d}{400} \frac{n}{10 d} = \frac{n}{4000},
\]
which is greater than $400 n / d$ for $d$ sufficiently large.  

\hfill $\Box$ 
\vspace{0.2cm}

\noindent
We next present the proof of Theorem \ref{t15}, which resembles that of Theorem
\ref{t11} but requires several additional ideas.
\vspace{0.2cm}

\noindent
{\bf Proof of Theorem \ref{t15}:}\,
 Let $c_1,c_2$ be the constants from Proposition \ref{p14} and put
\[
\delta = c_2/d.
\]

Let $\tau$ be a truth assignment consistent with $\rho$ and let $T = \tau^{-1}(\mathsf{t})$ and 
$B = \tau^{-1}(\mathsf{c})$. 
Let $H$ be the spanning  subgraph of $G$ in which an edge $(u,v)$
of $G$ is an edge of $H$ iff $\rho(u,v) = \mathsf{t}$. 
Let $V_1,V_2, \ldots ,V_s$ 
be the sets of vertices of the strongly connected  components 
(SCCs, for short) of 
$H$. 
\begin{clm}
\label{c25}
All the vertices of each $V_i$ are of the same type, that is, for
each $1 \leq i \leq s$, either $V_i \subset T$ or $V_i \subset B$.
\end{clm}

\noindent
{\bf Proof:}\, If $u \in T$ and $v$ is an out neighbor of
$u$ in $H$, then 
$v \in T$ (as $u$ reports so). If $v \in B$, and $u$ is an
in-neighbor of $v$ in $H$, then 
$u \in B$ (as $u$ reports that $u \in T$).  \hfill $\Box$

Call  an SCC  
of $H$ truthful if it is a subset of $T$, else it
is a subset of $B$ and we call it corrupt. 

Let $H'$ be the induced subgraph of $G$ on the set $T$ of all truthful 
vertices.
\begin{clm}
\label{c26}
Any SCC of $H'$ is also an SCC
of $H$.
\end{clm}

\noindent
{\bf Proof:}\, If $u,v \in T$  and $(u,v)$ is an edge of $G$, then
it is an edge of $H$ too.
Thus each SCC $C'$ of
$H'$ is contained in an SCC $C$ of $H$.  This SCC is
truthful, by Claim \ref{c25}, and cannot contain any additional 
truthful vertices as otherwise these belong to $C'$ as well. 
\hfill $\Box$

\begin{clm}
\label{c27}
The graph $H'$ contains an SCC of size at least
$|T|-2 \delta n > (1/2-2 \delta)n$. 
\end{clm}

\noindent
{\bf Proof:}\, Consider the component graph of $H'$: this is the 
directed graph $F$ whose vertices are all the SCCs of 
$H'$, where there is a directed 
edge  from $C$ to $C'$ iff there is some edge of  $H'$ from
some vertex of $C$ to some vertex of $C'$.  It is easy and well
known that this graph is a directed acyclic graph, and hence there
is a topological order of it, that is, a numbering
$C_1, C_2, \ldots ,C_r$ of the components  so that all edges
between different components  are of the form $(C_i,C_j)$ with
$i <j$. Order the vertices of $H'$ in a linear order according to
this topological order, where the vertices of $C_1$ come first
(in an arbitrary order),
those of $C_2$ afterwards, etc. Let $u_i$ be the vertex in place
$i$ according to this order $(1 \leq i \leq |T|)$. If the
vertices $u_{\delta n}$ and $u_{|T|-\delta n+1}$ 
belong to the same SCC, then this component is of size
at least $|T|-2 \delta n$ and we are done. Otherwise, the SCC
containing $u_{|T|/2}$ differs from either that containing
$u_{\delta n}$ or from that containing $u_{|T|-\delta n+1}$.
In the first case, the set $A$ of all SCCs up to that  containing
$u_{\delta n}$ is of size  at least $\delta n$, and the set
$B$ of all SCCs starting from that containing $u_{|T|/2}$ 
is of size at least $|T|/2 \geq n/4$. In addition there is no edge directed
from $B$ to $A$, contradicting the property of $G$. The second case
leads to a symmetric contradiction, establishing the claim.
\hfill $\Box$

Note that the above shows that if $|T|>(1/2+2\delta)n$ then
$H'$ and hence also $H$ must contain a SCC $C$ of size bigger than
$n/2$, which must be truthful.  Let $T'$ denote the set of vertices reachable from $C$ in $H$ and note that $T' \subset T$ and moreover $B' := N^{+}(T') \subset B$.
Let  $R = V \setminus (B' \cup T')$. 

By the expansion properties of $G$ it follows that 
$|R| \leq c_2 n / d$ and 
and that $|R| \leq c_3 |B'|/d \leq c_3 |B|/d$ 
for a constant $c_3=c_3(c_1,c_2)$. 

We next show that even if we only
assume that $|T|>n/2$ we can still identify correctly most of the
truthful vertices. 

By the last two claims it follows that $H$ contains at least one
SCC of size at least $(1/2-2 \delta)n \geq 3/8 n$.
If $H$ contains only one such component, then this component  
must consist of truthful agents, and we can identify all of them
(and hence also the types of all their out-neighbors). 
Otherwise, there is another SCC of size at least
$(1/2-\delta)n$, and as there is no room for more than two such
components, there are exactly two of them, say $V_1$ and $V_2$. Note
that by the properties of $G$ there are edges of $G$ from $V_1$
to $V_2$ and from $V_2$ to $V_1$, 
and hence it is impossible that both of them are
truthful components.  As one of them must be truthful, it follows
that exactly one  of them is truthful  and one is corrupt.
We next show that we can identify the types of both
components.

Recall that we have the SCCs of $H$, and the set $T$ of all truthful
vertices must be a union of a subset of these SCCs.
In addition, this set must be of size bigger than $n/2$ and must be
consistent with all reports of the vertices along every edge (in
the sense that for any edge $(u,v)$ with $u \in T$, the report of
$u$ on $v$ should be consistent with the actual type of $v$.)

\begin{clm}
\label{c28}
Given the strongly connected components $V_1,V_2, \ldots ,V_s$
of $H$ and the reports along each edge, either there is no 
union $I_1$ of SCCs including 
$V_1$ whose size exceeds $n/2$ so that $T=I_1,
B=V-I_1$ is
consistent with all reports along the edges, or there is no union
$I_2$ of SCCs including $V_2$ whose size exceeds $n/2$ so that 
$T=I_2,B=V-I_2$ is consistent with all reports along the edges. 
\end{clm}
\vspace{0.1cm}

\noindent
{\bf Proof:}\,  
Assume this is false, and let $I_1, I_2$ be as above.
By the above discussion we know that $I_1$ contains
$V_1$ but not $V_2$ and $I_2$ contains $V_2$ but not $V_1$.
Note that if some SCC $V_i$ is contained both in $I_1$ and in $I_2$ 
and there is any directed edge $(u,v)$ from $V_i$ to some other SCC
$V_j$, then if the report along this edge is that $v$ is truthful,
then $V_j$ must be  truthful  component in both $I_1$ and in $I_2$.
Similarly, if the report along this edge is $v \in B$, then $V_j$
must be outside $I_1$ and outside $I_2$. In particular, there are
no edges at all from $V_i$  to  $V_1$ or $V_2$ (as each of them
lies in exactly one of the two unions $I_1$, $I_2$). Let $J$ be the
set of all SCCs that are contained in both $I_1,I_2$. By the remark
above, for every edge $(u,v)$ from a vertex of $J$ to a vertex outside $J$,
the report along the edge must be  $v \in B$ (since otherwise $v$
would also be in an SCC which is truthful both in $I_1$ and in
$I_2$ and hence would be in $J$). Thus all edges $(u,v)$ as above
report $v \in B$, implying that all components outside $J$ 
to which there are
directed edges from vertices in $J$  belong to neither $I_1$ nor
$I_2$. By the properties of our graph the total size of these
components exceeds that of $J$, (as $|J| \leq 4 \delta n$ and all
out-neighbors of $J$ are outside $V_1,V_2$), and this shows that
the sum of the sizes of $I_1$ and $I_2$ is at  most
$$
2|J| +(|V|-|J|-|N^+(J)-J|) \leq |V|.
$$
Therefore  it cannot be that both $I_1$ and $I_2$ are of size bigger
than $|V|/2=n/2$, proving the claim. \hfill $\Box$

By the last claim it follows that one can identify the types of
the SCCs $V_1$ and $V_2$. This means that we can identify
at least $(1/2-2 \delta)n$ truthful vertices with no error. Recall
that this is
the case also when $H$ has only one SCC of size
at least $(1/2-\delta)n$. Having these truthful vertices, we also
know the types of all their out-neighbors. By the assumption on $G$
this gives the types of all vertices but at most $O(|B|/d)$,
completing the proof of the main part of the theorem.  
\vspace{0.1cm}

\noindent
The comment about the linear algorithm provided $|T|>(1/2+2\delta)n$
is clear. If we only assume that $|T|>|B|$ the proof
provides only a non-efficient algorithm  for deciding the types of
the SCCs $V_1$ and $V_2$. Indeed, we have to check
all $2^s$ possibilities of the types of each of the SCCs and see
which ones are consistent with all reports and are of total size
bigger than $n/2$.
By the proof above, only one of the two SCCs $V_1,V_2$ will appear 
among the truthful SCCs of such a possibility.
\hfill $\Box$

\section{Discussion and Open Problems}

The usefulness of expanders for corruption detection raises 
the natural question about the existence of good explicit 
spectral expanders with any desired number of nodes. 
After our work has been posted, explicit construction for 
such undirected expander graphs with any number of nodes were 
obtained in~\cite{A19}. 

It is
interesting to study in more detail the relation between expansion
and corruption detection.

\begin{quest}
Provide sharp criteria in terms of expansion and the fractional size of
the set $T$ for enabling corruption detection.
\end{quest}

For a weak result in the desired direction,  
consider the following argument.
We say that an undirected  graph $G$ is \emph{$\delta$-connected}
 if for every two disjoint sets
$A_1, A_2$ with $|A_1| \geq \delta n, | A_2 | \geq (1- 3 \delta) n$
there is at least one edge between $A_1$ and $A_2$. Note that the notion
of
$\delta$-connectedness   
is much weaker than expansion.  In particular a graph
$G$ can be
$\delta$-connected,   
yet at the same time have $\delta n/2 $
isolated vertices, while any $\delta$-good expander must be connected.

\begin{clm} \label{c:41} 
Suppose that $|T| = (1-\epsilon)n$ and the 
graph $G$ is $\epsilon$-connected. Then  
it is possible to identify 
$T' \subset T$ of size at least $(1-2\epsilon) n$. 
\end{clm}

\noindent
{\bf Proof:}\,
Let $E' \subset E$ be the set of edges both of 
whose end-points declare each other truthful. 
Recall that each connected components of $G'=(V,E')$ is either truthful or corrupt. 

Let $T_1, T_2,\ldots$ denote all the components of size at least $\epsilon
n$ in $G'$. Then we claim that if $T' = \cup T_i$ then $|T \setminus T'|
< \epsilon n$. Assume otherwise. Since all the connected components of
$T \setminus T'$ are of size at most $\epsilon n$, there exists $T''
\subset T \setminus T'$ of size in $[\epsilon n, 2 \epsilon n]$ with no
edges to 
$T \setminus T''$ whose size is in $[(1-3\epsilon)n ,(1-2\epsilon)
n]$. This is a contradiction to $\epsilon$-connectedness and the proof
follows.   \hfill $\Box$

To see that the conditions of Claim \ref{c:41} are tight up to constant
factors consider the star graph with $m$ leaves.  Assume that 
$|T| \leq
m-1$. Then it is easy to see that one cannot find even one member of $T$
if all vertices declare all their neighbors corrupt. On the other hand,
this example is (vacuously) $1/(4m)$ connected.  To get a non-trivial
example, one can replace each node with a complete graph $K_k$ and each
edge with a complete bipartite graph $K_{k,k}$ for an arbitrary $k$.

In a follow up work~\cite{JiMoRa:18}, the connection between expansion and corruption detection is formalized using the conjectured hardness of Small Set Expansion~\cite{RS}. Assuming the hardness of Small Set Expansion, 
it is shown that it is computationally hard to approximate the minimal number of nodes whose corruption makes identifying even one truthful node impossible. 

We conclude with a short discussion of a variant of the model. 
From the modeling perspective, it is interesting to consider 
probabilistic variants of the corruption detection problem. 
\begin{quest}
What is the effect of relaxing the assumption that truthful nodes
always report the status correctly? Suppose for example that each
truthful node reports the status of each of its neighbors independently
accurately with probability $1-\epsilon$.  Note that in this case it is
impossible to detect the status of an individual node with probability
one. However it is still desirable to find sets $T'$ and $B'$ such that
the symmetric difference $T \Delta T'$ and $B \Delta B'$ are small
with high probability. Under
what conditions can this be achieved? What are good algorithms for
finding $T'$ and $B'$?
\end{quest} 
See \cite{Alweiss:19} for results addressing this problem
obtained by Alweiss after our work has been posted.
\vspace{0.2cm}

\noindent
{\bf Acknowledgment:} We thank Gireeja Ranade for suggesting to consider
problems of corruption on networks, Peter Winkler for helpful
comments regarding the Byzantine Agreement problem, and two 
anonymous referees
for providing relevant references. We are also 
grateful to Laci Babai for 
numerous suggestions to improve the presentation as well as some of the results.

\bibliographystyle{plain}   
\bibliography{revref}

\end{document}